\documentclass[11pt]{amsart}
\usepackage[margin=1.25in]{geometry}
\usepackage{amsfonts,amsmath,amssymb,amsthm,graphics,subfig,enumerate,breqn,color, url}
\usepackage{enumitem}

\usepackage{hyperref}
\hypersetup{
    colorlinks=true,
    linkcolor=blue,
    filecolor=magenta,      
    urlcolor=cyan,
}
 
\makeatletter
\def\@tocline#1#2#3#4#5#6#7{\relax
  \ifnum #1>\c@tocdepth 
  \else
    \par \addpenalty\@secpenalty\addvspace{#2}%
    \begingroup \hyphenpenalty\@M
    \@ifempty{#4}{%
      \@tempdima\csname r@tocindent\number#1\endcsname\relax
    }{%
      \@tempdima#4\relax
    }%
    \parindent\z@ \leftskip#3\relax \advance\leftskip\@tempdima\relax
    \rightskip\@pnumwidth plus4em \parfillskip-\@pnumwidth
    #5\leavevmode\hskip-\@tempdima
      \ifcase #1
       \or\or \hskip 1em \or \hskip 2em \else \hskip 3em \fi%
      #6\nobreak\relax
    \hfill\hbox to\@pnumwidth{\@tocpagenum{#7}}\par
    \nobreak
    \endgroup
  \fi}
\makeatother

\newcommand{\Z}{\mathbb{Z}}

\newtheorem{lemma}{Lemma}[section]

\newtheorem{theorem}[lemma]{Theorem}
\newtheorem{corollary}[lemma]{Corollary}

\theoremstyle{definition}
\newtheorem{remark}{Remark}

\theoremstyle{remark}
\newtheorem*{proof*}{Proof}
\numberwithin{equation}{section}

\newenvironment{proofof}[1]{\noindent{\emph{Proof of {#1}.}}}{\qed\vspace{3ex}}

\newcommand\longvdots[1]{\raisebox{1em}{\rotatebox{-90}{\hbox to #1 {\dotfill}}}}
\long\def\symbolfootnote[#1]#2{\begingroup%
\def\thefootnote{\fnsymbol{footnote}}\footnote[#1]{#2}\endgroup}

\usepackage{cases}
\numberwithin{equation}{section}

\title{Counting subgroups of fixed order in finite abelian groups}
\author{Fikreab Admasu}
\address{Department of Mathematics, The Graduate Center, CUNY,  
365 Fifth Avenue, Room 4208, New York, NY 10016}
\email{fadmasu@gradcenter.cuny.edu}

\author{Amit Sehgal}
\address{Department of Mathematics, Government College, Birohar (Jhajjar) Haryana, India} 
\email{amit\_sehgal\_iit@yahoo.com}

\begin{document}
\maketitle

\begin{abstract}
We use recurrence relations to derive explicit formulas for counting the number of subgroups of given order (or index) in rank 3 finite abelian p-groups and use these to derive similar formulas in few cases for rank 4. As a consequence, we answer some questions by M. T$\ddot{a}$rn$\ddot{a}$uceanu in \cite{MT} and  L. T${\acute{o}}$th in \cite{LT}. We also  use other methods such as the method of fundamental group lattices introduced in \cite{MT} to derive a similar counting function in a special case of arbitrary rank finite abelian p-groups. 
\end{abstract}

\tableofcontents
\section{Introduction}

The subject of counting various kinds of subgroups of finite abelian groups has a long and rich history. For instance, in the early 1900's, Miller\cite{MI04} determined the number of cyclic subgroups of prime power order in a finite abelian p-group $G$, where p is a prime number. In about the same time, Hilton\cite{HH} found a necessary and sufficient condition for the existence of subgroups of a given order and then gave a general procedure for counting subgroups of order $p^r$ in a finite abelian p-group of order $p^n$, for positive integers $r\leq n.$ For recent work on counting subgroups of finite abelian p-groups, see \cite{GB96}, \cite{GBJW}, \cite[Sections 6 \& 7]{GLP}, \cite{HT}, \cite{AI}, \cite{MI39}, \cite{MT} and \cite{LT}. 

As an example, Hilton determined the total number of normal subgroups of index $p^2$ in any finite abelian p-group and noted that giving a general formula for the number of subgroups of order $p^r$ for every value of $r$ would be somewhat complicated. Progress has been made since then and there are now recursive formulas based on works of Stehling\cite{TS} and Birkhoff\cite{GB} that can compute the numbers for any $r.$ It can also be deduced using Hall's polynomials\cite{IGM} that the function that counts the number of subgroups is a polynomial in powers of $p.$ In addition, Shokuev\cite{VNS} finds an exact expression for the number of subgroups of any (possible) order of an arbitrary finite p-group. However, this requires subgroup structure of smaller rank subgroups in order to determine the expression. So far as we know, the explicit form of these counting polynomials has not been determined beyond that for rank 2 finite abelian p-groups. In \cite[Theorem 3]{GBJW01}, Bhowmik and Wu get a unimodality result on the coefficients and determine the leading coefficient of the polynomials in any rank. In this paper, we determine the explicit form of the counting polynomials in rank 3 and for some special values of $r$ in any rank. However, finding the general explicit polynomials in rank 4 or higher will be a tedious and difficult computation and we will illustrate this with examples.

By the fundamental theorem on the structure of finite abelian groups, any finite abelian group is isomorphic to a direct product of prime power order cyclic groups. Therefore, the number of subgroups in the finite abelian group is the product of the numbers of the subgroups contained in these cyclic groups\cite{RS}. This reduces the problem of counting subgroups in finite abelian groups to counting subgroups in finite abelian p-groups.

Let $\Z_{n}$ denote the cyclic group of order n. Given partitions $\lambda,\mu,\nu$ and $G$ an abelian $p$-group of type $\lambda  = (a_d,\dots,a_1),$ i.e., $G  \cong \Z/p^{a_1}\times \cdots \times \Z/p^{a_d},$ for some positive integers $a_1,\dots,a_d$ where $1\leq a_1\leq \dots \leq a_d,$ G has $g^\lambda_{\mu\nu}(p)$ subgroups $K$ such that $K$ has type $\mu$ and $G/K$ has type $\nu,$ where $g^\lambda_{\mu\nu}(X)\in \mathbb{Z}[X]$ is a Hall polynomial \cite{IGM}. Let 
\begin{equation}
  \label{no_of_subgrps}
 h_b^\lambda (p) = \sum_{|\nu| = b} \sum_{\mu}g^\lambda_{\mu\nu} (p)
\end{equation}
\noindent
where $\nu = (\nu_l,\dots,\nu_1)$ with $\nu_1\leq \nu_2\leq \dots \leq \nu_l$ and $|\nu| = \nu_1+\dots +\nu_l.$ 

Then $h_b^\lambda(p) = h_b^{(a_d,\dots,a_1)}(p)$ is the number of subgroups of index $p^b$ (or equivalently of order $p^{m-b}$ where $m=\sum_{i=1}^d a_i$) in the rank $d$ abelian $p$-group $\Z/p^{a_1}\times \cdots \times \Z/p^{a_d}.$ Moreover, $h_b^\lambda(p)$ is a polynomial as it is a sum of Hall polynomials. When the finite abelian $p$-group has rank 2, the number of subgroups of order $p^b$ in $\Z/p^{a_1}\times \Z/p^{a_2}$ is given by \cite[Theorem 3.3]{MT} as follows (for an alternate proof, see the end of Section \ref{sec:convolution}):
\begin{equation}
  \label{rank2ratfns}
  h_b^{(a_2,a_1)}(p)=\begin{cases}
  \frac{p^{b+1}-1}{p-1}  & \text{if } 0\le b \leq a_1\\    
  \frac{p^{a_1+1}-1}{p-1}  & \text{if } a_1\le b \leq a_2\\    
  \frac{p^{a_1+a_2-b+1}-1}{p-1}  & \text{if } a_2\le b \leq a_1+a_2.\\
\end{cases}
\end{equation}
\noindent
Except for rank 2 finite abelian p-groups and some other specific cases, there was no general direct closed formula (i.e. explicit polynomial depending on the type $\lambda$ and $b$; expressed above as a rational function in each case) that computes $h_b^\lambda(p)$ given $\lambda$ and $b.$ We determine in the next section the polynomials that give the number of subgroups of given order in rank 3 and few cases of rank 4 and one particular case for arbitrary rank finite abelian $p$-group. This answers part of a problem posed by M. T$\ddot{a}$rn$\ddot{a}$uceanu\cite[Problem 5.1]{MT}. As an application of the rational function formulas, we will prove a conjecture posed by Laszlo T${\acute{o}}$th in \cite{LT}. 
\section{Rank 3}
G. Birkhoff shows in \cite{GB} that $h_b^{(a_d,\dots,a_1)}$ is symmetric in $b$, that is, $h_b^{(a_d,\dots,a_1)} = h_{a_1+\dots+a_d-b}^{(a_d,\dots,a_1)}.$ So, for instance, in the three cases in equation (\ref{rank2ratfns}) for rank 2, the first case is symmetric with the third case while the second case is symmetric to itself. For example, when $d=2,$ the case $0\le b \leq a_1$ is symmetric to $a_2\le b \leq a_1+a_2$ because $$0\le b \leq a_1$$ 
$$-a_1\le -b \leq 0$$ 
$$a_2\leq a_1+a_2-b \leq a_1+a_2.$$
\noindent
Therefore, in order to determine the formula in the case $a_2\le b \leq a_1+a_2$ from the case when $0\le b \leq a_1$, we replace $b$ in the formula for $0\le b \leq a_1$ by $a_1+a_2-b.$ So, for instance in rank 2, it suffices to determine the compact form (as rational functions) of the polynomials in only two cases (first and second cases) and then use symmetry to deduce the remaining case (the last case). 

Let $\lambda=(a_3,a_2, a_1)$ and $\lambda^{'}=(a_2, a_1)$. The partition of the interval $[0,a_1+a_2+a_3]$ as in  Hironaka's paper \cite{YH} leads to 10 cases for rank 3. For each of these cases, there is a polynomial $h_b^{(a_1,a_2,a_3)}(p)$ (easily described as a rational function) enumerating the number of subgroups of index $p^b$ (or order $p^{a_1+a_2+a_3-b}$) in $\Z/p^{a_1}\times \Z/p^{a_2}\times \Z/p^{a_3}.$ 

There are at least four methods of deducing the formulas. The first method is to use a lemma \cite[Lemma 2.3]{YH} that Y. Hironaka deduces from a recurrence relation of T. Stehling in \cite{TS}. This is the method used in the proof of the theorem below as her lemma shows how to express the formulas in higher rank in terms of those of lower rank which we already know and this requires the least amount of computation. The second method was suggested to the first author by Gautam Chinta (the first author's doctoral advisor) while working on a related problem on subgroup growth as such subgroup counting formulas help compute smaller rank cotype zeta functions as in \cite{CKK}. This method involves using a series of contour integrals, the residue theorem and convolution of generating series to arrive at the formulas. The third method is based on an extension of results in the second author's previous paper in \cite{PKASSS}. This will be explicated another time. The last method, called the method of fundamental group lattices, was introduced by M. T$\ddot{a}$rn$\ddot{a}$uceanu in \cite{MT} and we have been able to apply it to derive a compact formula for counting subgroups in case 1 of all ranks. We demo the second method of proof in a special case at the end of the paper.

\begin{theorem}
\label{rank3}
For every $0\leq b \leq a_1+a_2+a_3$, the number $h_b^{(a_3,a_2,a_1)}(p)$ of all subgroups of order $p^b$ (equivalently  of index $p^{a_1+a_2+a_3-b}$) in the finite abelian $p$-group $\Z/p^{a_1}\times \Z/p^{a_2}\times \Z/p^{a_3}$ where $1\leq a_1\leq a_2 \leq a_3,$ is given by one of the following  polynomials expressed as rational functions:\\
Case 1: $0\le b \leq a_1$
\begin{equation*}
   h_b^{(a_3,a_2,a_1)}(p) = \frac{p^{2b+3} - p^{b+2} - p^{b+1} +1}{ (p-1)(p^2-1)}. 
\end{equation*}
Case 2:  $a_1\le b \leq a_2$ 
\begin{equation*}
h_b^{(a_3,a_2,a_1)}(p) = \frac{p^{b+a_1+3}+p^{b+a_1+2}-p^{2a_1+2}-p^{b+2}-p^{b+1}+1}{ (p-1)(p^2-1)}.
\end{equation*}
Case 3: $a_2 < b \le a_3 \leq a_1 + a_2$
\begin{equation*}
h_b^{(a_3,a_2,a_1)}(p) =\frac{(b-a_{2}+1)p^{a_2+a_1+3}+p^{a_2+a_{1}+2}-(b-a_2)p^{a_1+a_2+1}-p^{2a_1+2}-p^{b+2}-p^{b+1}+1}{(p-1)(p^{2}-1)}.
\end{equation*}
Case 4: $a_2 < b \le a_1 + a_2 \leq a_3$ 
$$h_b^{(a_3,a_2,a_1)}(p) = \frac{(b-a_{2}+1)p^{a_2+a_1+3}+p^{a_2+a_{1}+2}-(b-a_2)p^{a_1+a_2+1}-p^{2a_1+2}-p^{b+2}-p^{b+1}+1}{(p-1)(p^{2}-1)}.$$
Case 5: $a_1 + a_2 \le b \le a_3$ 
$$h_b^{(a_3,a_2,a_1)}(p) = \frac{(a_{1}+1)p^{a_2+a_1+3}+p^{a_2+a_{1}+2}-a_{1}p^{a_1+a_2+1}-p^{2a_1+2}-p^{a_{1}+a_{2}+2}-p^{a_{1}+a_{2}+1}+1}{(p-1)(p^{2}-1)}.$$
Case 6: $a_3 < b \leq a_{1}+a_{2} $
\begin{multline*}
  h_b^{(a_3,a_2,a_1)}(p) = \frac{(a_3-a_2+1)p^{a_2+a_1+3}+ 2p^{a_2+a_1+2}-(a_3-a_2-1)p^{a_2+a_1+1}-p^{a_3+a_2+a_1-b+2}}{(p-1)(p^2-1)}\\
  + \frac{-p^{a_3+a_2+a_1-b+1} -p^{2a_1+2}-p^{b+2}-p^{b+1}+1}{(p-1)(p^2-1)}.
\end{multline*}
The remaining cases are symmetric to one of the cases above, i.e., they can be obtained when we replace $b$ by $a_1+a_2+a_3-b.$\\
Case 7: $a_1 + a_2 < a_3  < b \le  a_1 + a_3$
\begin{multline*}
 h_b^{(a_3,a_2,a_1)}(p) = \frac{(a_1+a_3-b+1)p^{a_2+a_1+3}+p^{a_2+a_{1}+2}-(a_1+a_3-b)p^{a_1+a_2+1}-p^{2a_1+2}}{(p-1)(p^{2}-1)}\\ 
 + \frac{-p^{a_1+a_2+a_3-b+2}-p^{a_1+a_2+a_3-b+1}+1}{(p-1)(p^{2}-1)}.
\end{multline*}
Case 8: $ a_3  < a_1 + a_2 < b \le  a_1 + a_3$
\begin{multline*}
h_b^{(a_3,a_2,a_1)}(p) = \frac{(a_1+a_3-b+1)p^{a_2+a_1+3}+p^{a_2+a_{1}+2}-(a_1+a_3-b)p^{a_1+a_2+1}-p^{2a_1+2}}{(p-1)(p^{2}-1)}\\ + \frac{-p^{a_1+a_2+a_3-b+2}-p^{a_1+a_2+a_3-b+1}+1}{(p-1)(p^{2}-1)}.
\end{multline*}
Case 9: $a_1 + a_3 \le b \le  a_2 + a_3$
$$h_b^{(a_3,a_2,a_1)}(p) = \frac{p^{2a_1+a_2+a_3-b+3}+p^{2a_1+a_2+a_3-b+2}-p^{2a_1+2}-p^{a_1+a_2+a_3-b+2}-p^{a_1+a_2+a_3-b+1}+1}{(p-1)(p^{2}-1)}.$$
Case 10: $a_2 + a_3 \le b \le  a_1 + a_2 + a_3$
$$h_b^{(a_3,a_2,a_1)}(p) = \frac{p^{2a_1+2a_2+2a_3-2b+3}-p^{a_1+a_2+a_3-b+2}-p^{a_1+a_2+a_3-b+1}+1}{(p-1)(p^{2}-1)}.$$
\end{theorem}
\begin{proofof}{Theorem \ref{rank3}} To simplify notation, $h_b^{(a_3,a_2,a_1)}(p)$ will be denoted by $N_{b}(\lambda)$ as in the notation used in \cite{YH}. By \cite[Lemma 2.3]{YH}, we have the recursive formula 
\begin{equation}
    N_{b}(\lambda)=\sum_{i=0}^{b} p^{i}N_{i}(\lambda^{'}) - \sum_{|\lambda| +1-k}^{|\lambda^{'}|} p^{i}N_{i}(\lambda^{'}), 0\leq b\le |\lambda| 
\end{equation} 
\noindent where $\lambda = (a_d,\dots, a_1), \lambda^{'} = (a_{d-1},\dots, a_1), |\lambda| = a_d+\dots+a_1$ and the second summation appears only when $b > a_d.$ Thus, we compute each of the cases as follows: \\
Case 1: $0\le b \leq a_1$ \\
$N_{b}(\lambda)=\sum_{i=0}^{b} p^{i}N_{i}(\lambda^{'})=\sum_{i=0}^{b} p^{i} \frac{p^{i+1}-1}{p-1}=\sum_{i=0}^{b} \frac{p^{2i+1}-p^{i}}{p-1}=\frac{p^{2b+3}-p^{b+2}-p^{b+1}+1}{(p-1)(p^{2}-1)}.$\\
Case 2: $a_1\le b \leq a_2$ \\
$N_{b}(\lambda)=\sum_{i=0}^{a_1} p^{i}N_{i}(\lambda^{'})+\sum_{i=a_{1}+1}^{b} p^{i}N_{i}(\lambda^{'})=\frac{p^{2a_1+3}-p^{a_1+2}-p^{a_1+1}+1}{(p-1)(p^{2}-1)}+\sum_{i=a_{1}+1}^{b} p^{i} \frac{p^{a_{1}+1}-1}{p-1}$\\
$=\frac{p^{b+a_1+3}+p^{b+a_{1}+2}-p^{2a_1+2}-p^{b+2}-p^{b+1}+1}{(p-1)(p^{2}-1)}.$\\
Case 3: $a_2 < b \le a_3 \leq a_1 + a_2$\\ 
$N_{b}(\lambda)=\sum_{i=0}^{a_1} p^{i}N_{i}(\lambda^{'})+\sum_{i=a_{1}+1}^{a_{2}} p^{i}N_{i}(\lambda^{'})+\sum_{i=a_{2}+1}^{b} p^{i}N_{i}(\lambda^{'})$\\
$=\frac{p^{a_{2}+a_1+3}+p^{a_{2}+a_{1}+2}-p^{2a_1+2}-p^{a_{2}+2}-p^{a_{2}+1}+1}{(p-1)(p^{2}-1)}+\sum_{i=a_{2}+1}^{b} p^{i}\frac{p^{a_{1}+a_{2}+1-i}-1}{p-1}$\\
$=\frac{(b-a_{2}+1)p^{a_2+a_1+3}+p^{a_2+a_{1}+2}-(b-a_2)p^{a_1+a_2+1}-p^{2a_1+2}-p^{b+2}-p^{b+1}+1}{(p-1)(p^{2}-1)}.$\\
Case 4: $a_2 < b \le a_1 + a_2 \leq a_3$ \\
$N_{b}(\lambda)=\sum_{i=0}^{a_1} p^{i}N_{i}(\lambda^{'})+\sum_{i=a_{1}+1}^{a_{2}} p^{i}N_{i}(\lambda^{'})+\sum_{i=a_{2}+1}^{b} p^{i}N_{i}(\lambda^{'})$\\
$=\frac{p^{a_{2}+a_1+3}+p^{a_{2}+a_{1}+2}-p^{2a_1+2}-p^{a_{2}+2}-p^{a_{2}+1}+1}{(p-1)(p^{2}-1)}+\sum_{i=a_{2}+1}^{b} p^{i}\frac{p^{a_{1}+a_{2}+1-i}-1}{p-1}$\\
$=\frac{(b-a_{2}+1)p^{a_2+a_1+3}+p^{a_2+a_{1}+2}-(b-a_2)p^{a_1+a_2+1}-p^{2a_1+2}-p^{b+2}-p^{b+1}+1}{(p-1)(p^{2}-1)}.$\\
Case 5: $a_1 + a_2 \le b \le a_3$ \\
$N_{b}(\lambda)=\sum_{i=0}^{a_1} p^{i}N_{i}(\lambda^{'})+\sum_{i=a_{1}+1}^{a_{2}} p^{i}N_{i}(\lambda^{'})+\sum_{i=a_{2}+1}^{a_{1}+a_{2}} p^{i}N_{i}(\lambda^{'})+\sum_{i=a_{1}+a_{2}+1}^{b} p^{i}N_{i}(\lambda^{'})\\ - \sum_{i=a_{1}+a_{2}+a_{3}+1-b}^{a_{1}+a_{2}}p^{i}N_{i}(\lambda^{'})$\\
$=\frac{(a_{1}+1)p^{a_2+a_1+3}+p^{a_2+a_{1}+2}-a_{1}p^{a_1+a_2+1}-p^{2a_1+2}-p^{a_{1}+a_{2}+2}-p^{a_{1}+a_{2}+1}+1}{(p-1)(p^{2}-1)}.$\\
Case 6: $a_3 < b \leq a_{1}+a_{2} $ \\
$N_{b}(\lambda)=\sum_{i=0}^{a_1} p^{i}N_{i}(\lambda^{'})+\sum_{i=a_{1}+1}^{a_{2}} p^{i}N_{i}(\lambda^{'})+\sum_{i=a_{2}+1}^{a_{3}} p^{i}N_{i}(\lambda^{'})+\sum_{i=(+1}^{b} p^{i}N_{i}(\lambda^{'})\\ - \sum_{i=a_{1}+a_{2}+a3+1-b}^{a_{1}+a_{2}}p^{i}N_{i}(\lambda^{'})$\\
$=\frac{(a_{3}-a_{2}+1)p^{a_2+a_1+3}+p^{a_2+a_{1}+2}-(a_{3}-a_2)p^{a_1+a_2+1}-p^{2a_1+2}-p^{a_{3}+2}-p^{a_{3}+1}+1}{(p-1)(p^{2}-1)}+\sum_{i=a_{3}+1}^{b} p^{i}\frac{p^{a_{1}+a_{2}+1-i}-1}{p-1}$\\
$-\sum_{i=a_{1}+a_{2}+a_{3}+1-b}^{a_{1}+a_{2}}p^{i}\frac{p^{a_{1}+a_{2}+1-i}-1}{p-1}$\\
$=\frac{(a_3-a_2+1)p^{a_2+a_1+3}+ 2p^{a_2+a_1+2}-(a_3-a_2-1)p^{a_2+a_1+1}-p^{a_3+a_2+a_1-b+2}-p^{a_3+a_2+a_1-b+1}}{(p-1)(p^2-1)} \\
+ \frac{-p^{2a_1+2}-p^{b+2}-p^{b+1}+1}{(p-1)(p^2-1)}. $\\ 
The remaining cases are symmetric to one of the cases above, i.e., they can be obtained when we replace $b$ by $a_1+a_2+a_3-b.$\\
Case 7: $a_1 + a_2 < a_3  < b \le  a_1 + a_3$ Replace $b$ in case 4 by $a_1+a_2+a_3-b.$\\
$N_{b}(\lambda)= \frac{(a_1+a_3-b+1)p^{a_2+a_1+3}+p^{a_2+a_{1}+2}-(a_1+a_3-b)p^{a_1+a_2+1}-p^{2a_1+2}-p^{a_1+a_2+a_3-b+2}-p^{a_1+a_2+a_3-b+1}+1}{(p-1)(p^{2}-1)}.$\\
Case 8: $ a_3  < a_1 + a_2 < b \le  a_1 + a_3$ Replace $b$ in case 3 by $a_1+a_2+a_3-b.$\\
$N_{b}(\lambda) = \frac{(a_1+a_3-b+1)p^{a_2+a_1+3}+p^{a_2+a_{1}+2}-(a_1+a_3-b)p^{a_1+a_2+1}-p^{2a_1+2}-p^{a_1+a_2+a_3-b+2}-p^{a_1+a_2+a_3-b+1}+1}{(p-1)(p^{2}-1)}.$\\
Case 9: $a_1 + a_3 \le b \le  a_2 + a_3$ Replace $b$ in case 2 by $a_1+a_2+a_3-b.$\\
$N_{b}(\lambda) = \frac{p^{2a_1+a_2+a_3-b+3}+p^{2a_1+a_2+a_3-b+2}-p^{2a_1+2}-p^{a_1+a_2+a_3-b+2}-p^{a_1+a_2+a_3-b+1}+1}{(p-1)(p^{2}-1)}.$\\
Case 10: $a_2 + a_3 \le b \le  a_1 + a_2 + a_3$ Replace $b$ in case 1 by $a_1+a_2+a_3-b.$\\
$N_{b}(\lambda) = \frac{p^{2a_1+2a_2+2a_3-2b+3}-p^{a_1+a_2+a_3-b+2}-p^{a_1+a_2+a_3-b+1}+1}{(p-1)(p^{2}-1)}.$ 
\end{proofof}

\begin{remark} We can deduce all other cases using only case 6 as follows:\\
Case 1: replace $a_3,a_2$ and $a_1$ in  case 6 formula by $b.$\\
Case 2: replace $a_3$ and $a_2$ in  case 6 formula by $b.$\\
Cases 3 and 4: replace $a_3$ in  case 6 formula by $b.$\\
Case 5: replace $a_3$ and $b$  in  case 6 formula  by $a_2+a_1$ \\
Cases 7 and 8: replace $a_3$ in  case 6 formula by $a_1+a_2+a_3-b.$\\
Case 9: replace $a_3$ and $a_2$ in  case 6 formula by $a_1+a_2+a_3-b.$\\
Case 10: replace $a_3,a_2$ and $a_1$ in  case 6 formula by $a_1+a_2+a_3-b.$
\end{remark}
\section{Rank 4}
\subsection{Number of subgroups of given order}
While it was a somewhat lengthy but accomplishable task to determine the cases for rank 3, there does not appear to be a systematic way of determining all the cases for rank 4. We were able to list more than 22 cases some of which we list below. These formulas are also proved in exactly the same way as in Theorem \ref{rank3}. Here $\lambda=(a_4,a_3,a_2, a_1)$ and $\lambda^{'}=(a_3,a_2, a_1).$\\
When $0\le b \leq a_1,$ we get \\
$N_{b}(\lambda)=\sum_{i=0}^{b} p^{i}N_{i}(\lambda^{'})=\sum_{i=0}^{b} p^{i} \frac{p^{2i+3}-p^{i+2}-p^{i+1}+1}{(p-1)(p^{2}-1)}=\sum_{i=0}^{b}  \frac{p^{3i+3}-p^{2i+2}-p^{2i+1}+p^{i}}{(p-1)(p^{2}-1)}$\\
$=\frac{p^{3b+6}-p^{2b+5}-p^{2b+4}-p^{2b+3}+p^{b+3}+p^{b+2}+p^{b+1}-1}{(p-1)(p^{2}-1)(p^{3}-1)}.$\\
When $a_1\le b \leq a_2,$ we get 
$N_{b}(\lambda)=\sum_{i=0}^{a_1} p^{i}N_{i}(\lambda^{'})+\sum_{i=a_{1}+1}^{b} p^{i}N_{i}(\lambda^{'})$\\
$=\frac{p^{3a_{1}+6}-p^{2a_{1}+5}-p^{2a_{1}+4}-p^{2a_{1}+3}+p^{a_{1}+3}+p^{a_{1}+2}+p^{a_{1}+1}-1}{(p-1)(p^{2}-1)(p^{3}-1)}+\sum_{i=a_{1}+1}^{b} p^{i} \frac{p^{i+a_1+3}+p^{i+a_{1}+2}-p^{2a_1+2}-p^{i+2}-p^{i+1}+1}{(p-1)(p^{2}-1)}$\\
$=\frac{p^{2b+a_{1}+6}+p^{2b+a_{1}+5}+p^{2b+a_{1}+4}-p^{2a_{1}+b+5}-p^{2a_{1}+b+4}-p^{2a_{1}+b+3}+p^{3a_{1}+3}-p^{2b+5}-p^{2b+4}-p^{2b+3}+p^{b+3}+p^{b+2}+p^{b+1}-1}{(p-1)(p^{2}-1)(p^{3}-1)}.$\\
And as a final example, when $a_2\le b \leq min \{a_3, a_{1}+a_2 \},$ we have \\
$N_{b}(\lambda)=\sum_{i=0}^{a_1} p^{i}N_{i}(\lambda^{'})+\sum_{i=a_{1}+1}^{a_{2}} p^{i}N_{i}(\lambda^{'})+\sum_{i=a_{2}+1}^{b} p^{i}N_{i}(\lambda^{'})-\sum_{i=a_{1}+a_{2}+a_{3}+a_{4}+1-b}^{a_{1}+a_{2}+a_{3}}p^{i}N_{i}(\lambda^{'})$\\
$=\frac{p^{2a_{2}+a_{1}+6}+p^{2a_{2}+a_{1}+5}+p^{2a_{2}+a_{1}+4}-p^{2a_{1}+a_{2}+5}-p^{2a_{1}+a_{2}+4}-p^{2a_{1}+a_{2}+3}+p^{3a_{1}+3}}{(p-1)(p^{2}-1)(p^{3}-1)}\\
+\frac{-p^{2a_{2}+5}-p^{2a_{2}+4}-p^{2a_{2}+3}+p^{a_{2}+3}+p^{a_{2}+2}+p^{a_{2}+1}-1}{(p-1)(p^{2}-1)(p^{3}-1)}\\
+\sum_{i=a_{2}+1}^{b} p^{i} \frac{(i-a_{2}+1)p^{a_2+a_1+3}+p^{a_2+a_{1}+2}-(i-a_2)p^{a_1+a_2+1}-p^{2a_1+2}-p^{i+2}-p^{i+1}+1}{(p-1)(p^{2}-1)}$\\
$=\frac{(b+1-a_2)p^{a_{2}+a_{1}+b+6}+(b+1-a_{2})p^{a_{2}+a_{1}+b+5}+(a_{2}-b-1)p^{a_{2}+a_{1}+b+3}+(a_{2}-b-1)p^{a_{1}+a_{2}+b+2}}{(p-1)(p^{2}-1)(p^{3}-1)}\\
+\frac{p^{2a_{2}+a_{1}+4}+p^{a_{1}+2a_{2}+3}+p^{a_{1}+2a_{2}+2}}{(p-1)(p^{2}-1)(p^{3}-1)}$\\
$+\frac{-p^{2a_{1}+b+5}-p^{2a_{1}+b+4}-p^{2a_{1}+b+3}+p^{b+3}+p^{b+2}+p^{b+1}+p^{3a_{1}+3}-p^{2b+5}-p^{2b+4}-p^{2b+3}-1}{(p-1)(p^{2}-1)(p^{3}-1)}.$
\subsection{Total number of subgroups}
Let $\lambda = (a_4,a_3,a_2,a_1)$ and let $N(\lambda)$ denote the total number of subgroups of $\Z/p^{a_1}\times \Z/p^{a_2}\times \Z/p^{a_3}\times \Z/p^{a_4}$ where $1\le a_1\le a_2\le a_3\le a_4,$ and $n= a_1+a_2+a_3+a_4.$ Then 
\begin{equation*}
N(\lambda) = \sum_{b=0}^n N_b(\lambda).
\end{equation*}

A conjecture of L. T$\acute{o}$th in \cite[Conjecutre 10]{LT} claims that the degree of the polynomial $N(m,m,m,m)$, i.e. when $\lambda = (m,m,m,m),$ is $4m$ and that it has leading coefficient of 1. We prove that this follows from our formulas in Theorem \ref{rank3}. First, let's state a corollary for the special case of a rank 3 abelian p-group of type $(m,m,m).$ The proof of the corollary follows immediately from setting $a_1=a_2=a_3=m$ in the theorem. 

\begin{corollary}\label{mmm}
For every $0\leq b \leq 3m$, the number $h_b^{(m,m,m)}(p)$ of all subgroups of order $p^b$ in the finite abelian $p$-group $\Z/p^{m}\times \Z/p^{m}\times \Z/p^{m}$ where $1\leq m,$ is given by one of the following polynomials expressed as rational functions:\\
Case 1: $0\le b \leq m$
$$h_b^{(m,m,m)}(p) = \frac{p^{2b+3} - p^{b+2} - p^{b+1} +1}{(p-1)(p^2-1)}.$$
Case 2:  $m \le b \leq 2m$
$$h_b^{(m,m,m)}(p) = \frac{p^{2m+3}+p^{2m+2}+p^{2m+1}-p^{3m+2-b}-p^{3m+1-b}-p^{b+2}-p^{b+1}+1}{ (p-1)(p^2-1)}.$$
Case 3: $2m \le b \le 3m$ 
$$h_b^{(m,m,m)}(p) =\frac{p^{6m+3-2b}-p^{3m+2-b}-p^{3m+1-b}+1}{(p-1)(p^{2}-1)}.$$
\end{corollary}
\begin{theorem}\label{Conj10}
The the total number of subgroups $N(m,m,m,m)$ of the rank 4 finite abelian p-group $\Z/p^{m}\times \Z/p^{m}\times \Z/p^{m}\times \Z/p^{m}$ where $1\le m$ is given by a polynomial expressed as a rational function as below:
\begin{multline*}
 N(m,m,m,m) = \sum_{b=1}^{4m} N_b(m,m,m,m)\\
 = \frac{ {\left(p^{2} + p + 1\right)}^{3} {\left(p^{2} + 1\right)} p^{4 \, m + 2} - {\left({\left(2 \, m + 3\right)} p^{3} - 2 \, m - 1\right)} {\left(p^{3} + p\right)} {\left(p + 1\right)}^{3} p^{3 \, m} }{{\left(p^{2} - 1\right)}^{2} {\left(p^3 - 1\right)}^{2}}\\
 - \frac{4 \, m p^{4} + 4 \, m p^{3} + 7 \, p^{4} + 9 \, p^{3} - 4 \, m p + 6 \, p^{2} - 4 \, m + p - 1}{{\left(p^{2} - 1\right)}^{2} {\left(p^3 - 1\right)}^{2}}.
\end{multline*}
In particular, the degree of the polynomial is $4m$ and its leading coefficient is 1.

\end{theorem}

\begin{proofof}{Theorem \ref{Conj10}}
We break the interval $[0,4m]$ into 4 sub-intervals and use \cite[Lemma 3.2]{YH} and the corollary above in each case to derive a polynomial expression for the number of subgroups of order $p^b$ in a rank 4 abelian p-group of type $(m,m,m,m)$. Summing over these polynomials will result in the required polynomial expression.\\
Case 1: $0\le b\le m$
\begin{align*}
N_b(m,m,m,m) &= \sum_{i=0}^{b} p^{i}N_{i}(m,m,m)=\sum_{i=0}^{b} p^{i} \frac{p^{2i+3} - p^{i+2} - p^{i+1} +1}{(p-1)(p^2-1)}\\
&= \frac{p^{3 \, b + 6} - p^{2 \, b + 5} - p^{2 \, b + 4} - p^{2 \, b + 3} + p^{b + 3} + p^{b + 2} + p^{b + 1} - 1}{{\left(p - 1\right)} {\left(p^2 - 1\right)} {\left(p^3 - 1\right)}}. 
\end{align*}
Case 2: $m < b\le 2m$
\begin{align*}
N_b(m,m,m,m) &= \sum_{i=0}^{m} p^{i}N_{i}(m,m,m) + \sum_{i=m+1}^{b} p^{i}N_{i}(m,m,m) - \sum_{i=4m+1-b}^{3m} p^{i}N_{i}(m,m,m)\\
&=  \sum_{i=0}^{m} p^{i}\frac{p^{2i+3} - p^{i+2} - p^{i+1} +1}{(p-1)(p^2-1)}\\ &+ \sum_{i=m+1}^{b} p^{i}\frac{p^{2m+3}+p^{2m+2}+p^{2m+1}-p^{3m+2-i}-p^{3m+1-i}-p^{i+2}-p^{i+1}+1}{ (p-1)(p^2-1)}\\
&- \sum_{i=4m+1-b}^{3m} p^{i}\frac{p^{6m+3-2i}-p^{3m+2-i}-p^{3m+1-i}+1}{(p-1)(p^{2}-1)}\\
&= -\frac{p^{2 \, b + 5} + p^{2 \, b + 4} + p^{2 \, b + 3} - p^{b + 2 \, m + 6} - p^{b + 2 \, m + 5} - 2 \, p^{b + 2 \, m + 4}}{(p-1)(p^2-1)(p^3-1)} \\
&+\frac{ - p^{b + 2 \, m + 3} - p^{b + 2 \, m + 2} - p^{b + 3} - p^{b + 2}}{(p-1)(p^2-1)(p^3-1)} \\
&-\frac{ - p^{b + 1} - p^{-b + 4 \, m + 3} - p^{-b + 4 \, m + 2} - p^{-b + 4 \, m + 1} + p^{3 \, m + 5} + 2 \, p^{3 \, m + 4}}{(p-1)(p^2-1)(p^3-1)}\\ 
&+ \frac{2 \, p^{3 \, m + 3} + 2 \, p^{3 \, m + 2} + p^{3 \, m + 1} + 1}{(p-1)(p^2-1)(p^3-1)}.
\end{align*}
Case 3: $2m < b \le 3m$
\begin{align*}
N_b(m,m,m,m) &= \sum_{i=0}^{m} p^{i}N_{i}(m,m,m) + \sum_{i=m+1}^{2m} p^{i}N_{i}(m,m,m) + \sum_{i=2m+1}^{b} p^{i}N_{i}(m,m,m)\\ &- \sum_{i=4m+1-b}^{2m} p^{i}N_{i}(m,m,m) -\sum_{i=2m+1}^{3m} p^{i}N_{i}(m,m,m)\\
&=  \sum_{i=0}^{m} p^{i}\frac{p^{2i+3} - p^{i+2} - p^{i+1} +1}{(p-1)(p^2-1)}\\
&+ \sum_{i=m+1}^{2m} p^{i}\frac{p^{2m+3}+p^{2m+2}+p^{2m+1}-p^{3m+2-i}-p^{3m+1-i}-p^{i+2}-p^{i+1}+1}{ (p-1)(p^2-1)}\\
&+ \sum_{i=2m+1}^{b} p^{i}\frac{p^{6m+3-2i}-p^{3m+2-i}-p^{3m+1-i}+1}{(p-1)(p^{2}-1)}\\
&- \sum_{i=4m+1-b}^{2m} p^{i}\frac{p^{2m+3}+p^{2m+2}+p^{2m+1}-p^{3m+2-i}-p^{3m+1-i}-p^{i+2}-p^{i+1}+1}{ (p-1)(p^2-1)}\\ &- \sum_{i=2m+1}^{3m} p^{i}\frac{p^{6m+3-2i}-p^{3m+2-i}-p^{3m+1-i}+1}{(p-1)(p^{2}-1)}\\
&=\frac{p^{b + 3} + p^{b + 2} + p^{b + 1} + p^{-b + 6 \, m + 6} + p^{-b + 6 \, m + 5} + 2 \, p^{-b + 6 \, m + 4}}{(p-1)(p^2-1)(p^3-1)}\\ &+ \frac{p^{-b + 6 \, m + 3} + p^{-b + 6 \, m + 2} + p^{-b + 4 \, m + 3}}{(p-1)(p^2-1)(p^3-1)}\\
&+ \frac{p^{-b + 4 \, m + 2} + p^{-b + 4 \, m + 1} - p^{-2 \, b + 8 \, m + 5} - p^{-2 \, b + 8 \, m + 4} - p^{-2 \, b + 8 \, m + 3} - p^{3 \, m + 5} - 2 \, p^{3 \, m + 4}}{(p-1)(p^2-1)(p^3-1)}\\
&+ \frac{- 2 \, p^{3 \, m + 3} - 2 \, p^{3 \, m + 2} - p^{3 \, m + 1} - 1}{(p-1)(p^2-1)(p^3-1)}.
\end{align*}
Case 4: $3m < b\le 4m$
\begin{align*}
N_b(m,m,m,m) &= \sum_{i=0}^{m} p^{i}N_{i}(m,m,m) + \sum_{i=m+1}^{2m} p^{i}N_{i}(m,m,m) + \sum_{i=2m+1}^{3m} p^{i}N_{i}(m,m,m) \\
&- \sum_{i=4m+1-b}^{m} p^{i}N_{i}(m,m,m) -\sum_{i=m+1}^{2m} p^{i}N_{i}(m,m,m) - \sum_{i=2m+1}^{3m} p^{i}N_{i}(m,m,m)\\
&=  \sum_{i=0}^{m} p^{i}\frac{p^{2i+3} - p^{i+2} - p^{i+1} +1}{(p-1)(p^2-1)}\\
&+ \sum_{i=m+1}^{b} p^{i}\frac{p^{2m+3}+p^{2m+2}+p^{2m+1}-p^{3m+2-i}-p^{3m+1-i}-p^{i+2}-p^{i+1}+1}{ (p-1)(p^2-1)}
\end{align*}
\begin{align*}
&+ \sum_{i=2m+1}^{3m} p^{i}\frac{p^{6m+3-2i}-p^{3m+2-i}-p^{3m+1-i}+1}{(p-1)(p^{2}-1)}\\
&- \sum_{i=4m+1-b}^{m} p^{i}\frac{p^{2i+3} - p^{i+2} - p^{i+1} +1}{(p-1)(p^2-1)}\\
&- \sum_{i=m+1}^{2m} p^{i}\frac{p^{2m+3}+p^{2m+2}+p^{2m+1}-p^{3m+2-i}-p^{3m+1-i}-p^{i+2}-p^{i+1}+1}{ (p-1)(p^2-1)}\\ 
&- \sum_{i=2m+1}^{3m} p^{i}\frac{p^{6m+3-2i}-p^{3m+2-i}-p^{3m+1-i}+1}{(p-1)(p^{2}-1)}\\
&=\frac{p^{-b + 4 \, m + 3} + p^{-b + 4 \, m + 2} + p^{-b + 4 \, m + 1} - p^{-2 \, b + 8 \, m + 5}}{(p-1)(p^2-1)(p^3-1)}\\
&+ \frac{- p^{-2 \, b + 8 \, m + 4} - p^{-2 \, b + 8 \, m + 3} + p^{-3 \, b + 12 \, m + 6} - 1}{(p-1)(p^2-1)(p^3-1)}.
\end{align*}
Now summing over the 4 cases, we get that
\begin{align*}
N(m,m,m,m)&=\sum_{b=0}^{4m} N_b(m,m,m,m)\\
&= \frac{ {\left(p^{2} + p + 1\right)}^{3} {\left(p^{2} + 1\right)} p^{4 \, m + 2} - {\left({\left(2 \, m + 3\right)} p^{3} - 2 \, m - 1\right)} {\left(p^{3} + p\right)} {\left(p + 1\right)}^{3} p^{3 \, m} }{{\left(p^{2} - 1\right)}^{2} {\left(p^3 - 1\right)}^{2}}\\
&- \frac{4 \, m p^{4} + 4 \, m p^{3} + 7 \, p^{4} + 9 \, p^{3} - 4 \, m p + 6 \, p^{2} - 4 \, m + p - 1}{{\left(p^{2} - 1\right)}^{2} {\left(p^3 - 1\right)}^{2}}.
\end{align*}
Now, we can see that the degree of the polynomial is $4m$ and the leading coefficient is 1, as conjectured by L. T$\acute{o}$th.
\end{proofof}
\begin{remark} In an unpublished paper \cite{CCL}, Chew, Chin, and Lim derive an explicit formula for the number of subgroups of a finite abelian p-group of rank 4. Here we will use the formula to reprove the above theorem and also prove one more conjecture by T$\acute{o}$th in \cite[Conjecture 9]{LT}. 
\begin{theorem} (Chew, Chin, and Lim) Let $1 \le w \le x \le y \le z$. The number of subgroups of $\Z/p^{w}\times \Z/p^{x}\times \Z/p^{y}\times \Z/p^{z}$ is
\begin{align*}
N(w,x,y,z) &= \sum_{i=0}^{w-1}\sum_{j=0}^{i} [(w+x+y+z-4i+1)(2i-2j+1)p^{3i+j}\\
&+(w+x+y+z-4i-1)(2i-2j+1)p^{3i+j+1}\\
&+ 2(w+x+y+z-4i-2)(i-j+1)p^{3i+j+2}]\\
&+\sum_{i=0}^{w}\sum_{j=w}^{x-1} (w+j-2i+1)[(x+y+z-3j+1)p^{w+2j+i}\\
\end{align*}
\begin{align*}
&+(x+y+z-3j-1)p^{w+2j+i+1}] \\
&+\sum_{i=0}^{w}\sum_{j=x}^{y} (y+z-2j+1)(w+x+-2i+1)p^{w+x+i+j}.
\end{align*}
\end{theorem}
Since all coefficients are positive, none of the terms will cancel out and the degree of the polynomial in $p$ can be easily determined. By comparing the exponents, we can see that the highest degree appears in the last double sum. Setting $i=w$ and $j=y,$ we get that the leading term of $N(w,x,y,z)$ is $(z-y+1)(x-w+1)p^{2w+x+y},$ with degree $2w+x+y,$ confirming conjecture 9 in \cite{LT}. 

Moreover, setting $w=x=y=z=m,$ in the above leading term, we get that the degree of the leading term in $N(m,m,m,m)$ is $4m,$ confirming Conjecture 10 in the same paper. 
\end{remark}

\begin{remark} What is new in Theorem \ref{Conj10} is the explicit polynomial. Else, both conjectures 9 and 10 of Toth follow from Theorem 3 in \cite{GBJW01} as they have determined the leading coefficient and degree of the counting polynomial for all ranks.
\end{remark}

\section{A special case of any rank}
The recursive counting of the subgroups used in the proof of Theorem \ref{rank3} based on \cite[Lemma 2.3]{YH} enabled us to deduce a compact rational form for the polynomials in all the cases for ranks 2 and 3. Here we will demo the method of fundamental group lattices in one case for all ranks. We remark that it is difficult to use this method to give an explicit formula in the remaining cases (those that are not symmetric to Case 1) of rank 3 or higher.

Below we extend T$\ddot{a}$rn$\ddot{a}$uceanu's result for counting subgroups as in \cite[Theorem 3.3]{MT} for case 1 to any rank. 
\begin{theorem} 
\label{funlatt} Let $0\leq b \leq a_1.$ Then 
\begin{equation}
h_b^{(a_k,\dots,a_1)}(p) = \prod_{i=2}^k\frac{p^{b+i-1}-1}{p^{i-1}-1}. 
\end{equation} Similarly, if $a_2+\dots+a_k\le b \leq a_1+\dots+a_k$, then replace $b$ in the above formula by $a_1+\dots+a_k-b$ to get 
\begin{equation}
h_b^{(a_k,\dots,a_1)}(p) = \prod_{i=2}^k\frac{p^{a_1+\dots+a_k-b+i-1}-1}{p^{i-1}-1}.  \end{equation}
\end{theorem}
\begin{proofof}{Theorem \ref{funlatt}} Following the notation in \cite{MT}, let $A= (a_{ij})$ be a solution of $(\ast)$ below corresponding to a subgroup of order $p^{a_1+\dots+a_k-b}$ in $\Z/p^{a_1}\times \Z/p^{a_2}\times \cdots \times \Z/p^{a_k}.$ \\

$$(\ast) \begin{cases} \text{i) } a_{ij} = 0 \text{ for any } i>j\\
\text{ii) } 0\leq a_{1j},a_{2j},\dots,a_{j-1,j} < a_{jj} \text{ for any } j \in \{1,2,\dots, k\}\\
\text{iii) } \text{1) } a_{11}|p^{a_1}\\
\hspace{0.2in}\text{ 2) } a_{22}|(p^{a_2},p^{a_1}\frac{a_{12}}{a_{11}})\\
\hspace{0.2in}\text{ 3) } a_{33}|(p^{a_3},p^{a_2}\frac{a_{23}}{a_{22}},p^{a_1}\frac{\begin{vmatrix}
a_{12} &a_{13} \\ 
a_{22} & a_{23}
\end{vmatrix}}{a_{22}a_{11}})\\
\vspace{0.1in}\\
\hspace{1in} \longvdots{1.5em}\\
\hspace{0.2in}\text{ k) } a_{kk}|(p^{a_k},p^{a_{k-1}}\frac{a_{k-1,k}}{a_{k-1,k-1}},p^{a_{k-2}}\frac{\begin{vmatrix}
a_{k-2,k-1} &a_{k-2,k} \\ 
a_{k-1,k-1} & a_{k-1,k}
\end{vmatrix}}{a_{k-1,k-1}a_{k-2,k-2}},\dots,p^{a_1}\frac{\begin{vmatrix}
a_{12} & a_{13} & \cdots & a_{1,k} \\ 
a_{22} & a_{23} & \cdots & a_{2,k} \\ 
\textbf{.}  & \textbf{.}  & & \textbf{.}  \\ 
\textbf{.}  & \textbf{.}  & & \textbf{.}  \\ 
\textbf{.}  & \textbf{.}  & & \textbf{.}  \\ 
a_{k-1,2} & a_{k-1,3} & \cdots & a_{k-1,k} 
\end{vmatrix}}{a_{k-1,k-1}a_{k-2,k-2}\dots a_{11}}).
\end{cases} $$

Put $a_{11}=p^{i_1}$ where $0\leq i_1\leq a_1.$ By the first remark on pp. 376 of \cite{MT}, the order of a subgroup corresponding to $A= (a_{ij})$ is $$\frac{p^{\sum_{i=1}^k a_i}}{\prod_{i=1}^k a_{ii}} = p^{a_1+\dots+a_k-b}$$
Therefore, $a_{11}\cdots a_{kk}=p^b.$ Since $a_{11}=p^{i_1},$ we have $a_{22}\cdots a_{kk}=p^{b-i_1}.$ Let\\

\hspace{2in} $a_{22} = p^{i_2}$ where $0\leq i_2 \leq  b-i_1,$ \\ 

\hspace{2in} $a_{33}=p^{i_3}$ where $0\leq i_3\leq b-i_1-i_2,$\\

\hspace{2in} $\dots $ \\

\hspace{2in} $a_{k-1,k-1}=p^{i_{k-1}}$ where $0\leq i_{k-1} \leq b-i_1-\dots - i_{k-2}.$ \\
Then $a_{kk}=p^{b-i_1-\dots - i_{k-1}}.$ Now, the conditions become
$$p^{i_2}|(p^{a_2},p^{a_1-i_1}a_{12}) = p^{a_1-i_1}(p^{a_2-a_1+i_1},a_{12})$$
$$p^{i_3}|(p^{a_3},p^{a_2-i_2}a_{12}, p^{a_1-i_1-i_2}(a_{12}a_{23}-p^{i_2}a_{13})$$
$$\dots$$
$$p^{b-i_1-\dots - i_{k-1}}|(p^{a_k},p^{a_{k-1}}\frac{a_{k-1,k}}{a_{k-1,k-1}},p^{a_{k-2}}\frac{\begin{vmatrix}
a_{k-2,k-1} &a_{k-2,k} \\ 
a_{k-1,k-1} & a_{k-1,k}
\end{vmatrix}}{a_{k-1,k-1}a_{k-2,k-2}},\dots,p^{a_1}\frac{\begin{vmatrix}
a_{12} & a_{13} & \cdots & a_{1,k} \\ 
a_{22} & a_{23} & \cdots & a_{2,k} \\ 
\textbf{.}  & \textbf{.}  & & \textbf{.}  \\ 
\textbf{.}  & \textbf{.}  & & \textbf{.}  \\ 
\textbf{.}  & \textbf{.}  & & \textbf{.}  \\ 
a_{k-1,2} & a_{k-1,3} & \cdots & a_{k-1,k} 
\end{vmatrix}}{a_{k-1,k-1}a_{k-2,k-2}\dots a_{11}}).$$
These conditions are satisfied by all
$$a_{12}< a_{22}=p^{i_2},$$
$$a_{13},a_{23} < a_{33}=p^{i_3},$$
$$\cdots$$
$$a_{1,k-1},\dots,a_{k-2,k-1}< a_{k-1,k-1}=p^{i_{k-1}}$$
$$a_{1,k},\dots,a_{k-1,k}< a_{k,k}=p^{b-i_1-\cdots-i_{k-1}}.$$
So, we get $$p^{i_2}\cdot (p^{i_3})^2 \cdots (p^{i_{k-1}})^{k-2}\cdot (p^{b-i_1-\cdots-i_{k-1}})^{k-1} = p^{(k-1)b-(k-1)i_1-(k-2)i_2-\cdots -2i_{k-2}-i_{k-1}}$$ distinct solutions of $(\ast)$ for given $a_{11}=p^{i_1},a_{22}=p^{i_2},\cdots,a_{k-1,k-1}=p^{i_{k-1}},$ and  $a_{k,k}=p^{b-i_1-\cdots-i_{k-1}}.$ The number of subgroups of index $p^b$ (or order $p^{a_1+\dots+a_k-b}$) in $\Z/p^{a_1}\times \cdots \times \Z/p^{a_k}$ is then 
\begin{equation}
h_b^{(a_k,\dots,a_1)}(p) = \sum_{i_1=0}^{b}\sum_{i_2=0}^{b-i_1}\cdots \sum_{i_{k-1}=0}^{b-i_1-\cdots-i_{k-2}} p^{(k-1)b-(k-1)i_1-(k-2)i_2-\cdots -2i_{k-2}-i_{k-1}}.
\end{equation}

Now we will use induction to prove that 
\begin{equation}\label{anyrankclaim}
\sum_{i_1=0}^{b}\sum_{i_2=0}^{b-i_1}\cdots \sum_{i_{k-1}=0}^{b-i_1-\cdots-i_{k-2}} p^{(k-1)b-(k-1)i_1-(k-2)i_2-\cdots -2i_{k-2}-i_{k-1}} = \prod_{i=2}^k\frac{p^{b+i-1}-1}{p^{i-1}-1}. 
\end{equation}

By Theorem \ref{rank3} and previous discussion, equation (\ref{anyrankclaim}) is true for when $k=2,3.$ Let us assume that it is true for rank k. Now consider the left hand side of equation (\ref{anyrankclaim})  for rank $k+1,$ that is, 
$$\sum_{i_1=0}^{b}\sum_{i_2=0}^{b-i_1}\cdots \sum_{i_{k-1}=0}^{b-i_1-\cdots-i_{k-2}} \sum_{i_{k}=0}^{b-i_1-\cdots-i_{k-2}-i_{k-1}}p^{(k)b-(k)i_1-(k-1)i_2-\cdots -3i_{k-2}-2i_{k-1}-i_{k}}  $$\\
and put $\alpha=b-i_1$. Then \\
$$\sum_{i_1=0}^{b}p^{b-i_1}\sum_{i_2=0}^{\alpha}\cdots \sum_{i_{k-1}=0}^{\alpha-\cdots-i_{k-2}} \sum_{i_{k}=0}^{\alpha-\cdots-i_{k-2}-i_{k-1}}p^{(k-1)\alpha-(k-1)i_2-\cdots -3i_{k-2}-2i_{k-1}-i_{k} }$$
\begin{align*}
&= \sum_{i_1=0}^{b}p^{b-i_1} \prod_{i=3}^{k+1}\frac{p^{b+i-1}-1}{p^{i-1}-1}\\
&= \prod_{i=3}^{k+1}\frac{p^{b+i-1}-1}{p^{i-1}-1}\sum_{i_1=0}^{b}p^{b-i_1}\\
&= \prod_{i=2}^{k+1}\frac{p^{b+i-1}-1}{p^{i-1}-1}.
\end{align*} This completes the proof. 
\end{proofof}
\begin{remark} An application of the method used in Theorem \ref{rank3} also shows that we can rederive the above result in a simpler way as follows. Let $\lambda=(a_k,\dots, a_1), \lambda^{'}=(a_{k-1},\dots, a_1),$ and  $0\le b \leq a_1.$ Then a repeated application of the recurrence (\ref{Stehrec}) for $N_b(\lambda)$ shows that  
\begin{align*}
N_{b}(\lambda) &= \sum_{i_1=0}^{b} p^{i_1}N_{i_1}(\lambda^{'})\\
 &=\sum_{i_1=0}^{b} p^{i_1}\sum_{i_2=0}^{i_1} p^{i_2}N_{i_2}(\lambda^{''})\\
 &=\dots \\
 &=\sum_{i_1=0}^{b} \sum_{i_2=0}^{i_1} \dots\sum_{i_{k-1}=0}^{i_{k-2}} p^{i_1+\dots +i_{k-1}} N_{i_{k-1}}(a_1)\\
 &=\sum_{i_1=0}^{b} \sum_{i_2=0}^{i_1} \dots\sum_{i_{k-1}=0}^{i_{k-2}} p^{i_1+\dots+i_{k-1}}.
\end{align*}

Now let $j_t = b - \sum_{l=1}^t i_t$ for $1\le t\le k-2.$ Then $i_t = j_{t-1} - j_t$ for $1\le t\le k-2,$ letting $j_0=b.$ As $i_t$ varies between 0 and $j_{t-1},$ $j_t$ will also vary between 0 and $j_{t-1}.$ Also $(k-1)b-(k-1)i_1-(k-2)i_2-\cdots -2i_{k-2}-i_{k-1}= j_1+\dots+j_{k-1}.$ Hence,
\begin{align*}
\sum_{j_1=0}^{b} \sum_{j_2=0}^{j_1} \dots\sum_{j_{k-1}=0}^{j_{k-2}} p^{j_1+\dots+j_{k-1}} 
&= \sum_{i_1=0}^{b} \sum_{i_2=0}^{j_1} \dots\sum_{i_{k-1}=0}^{j_{k-2}} p^{j_1+\dots+j_{k-1}} \\
&=\sum_{i_1=0}^{b}\sum_{i_2=0}^{b-i_1}\cdots \sum_{i_{k-1}=0}^{b-i_1-\cdots-i_{k-2}} p^{(k-1)b-(k-1)i_1-(k-2)i_2-\cdots -2i_{k-2}-i_{k-1}},
\end{align*} which is the same as the left hand side of equation (\ref{anyrankclaim}). Therefore, 
\begin{equation}
\sum_{i_1=0}^{b} \sum_{i_2=0}^{i_1} \dots\sum_{i_{k-1}=0}^{i_{k-2}} p^{i_1+\dots+i_{k-1}} =  \prod_{i=2}^k\frac{p^{b+i-1}-1}{p^{i-1}-1}.   
\end{equation}

\end{remark}
\section{The method of convolution} \label{sec:convolution}
Here, we will use convolution of generating series and a recurrence relation of T. Stehling\cite{TS} to prove the rational function formula in case 1 of rank 2. This is only to demo the method and the rest of the cases can be handled similarly. To ease the use of the recurrence relation, let's adopt Stehling's notation for type $\alpha = (\alpha_1,\dots,\alpha_d)$ where $ \alpha_1\ge \cdots \alpha_d\ge 0$ which reverses the order in the notation we used in previous sections. 

Let $0\le b \leq \alpha_2$ (i.e. case 1 of rank 2), then
\begin{equation}
\label{c1rk2}
h_b^{(\alpha_1,\alpha_2)}(p) = \frac{ p^{b+1} -1}{p-1}. 
\end{equation}
The symmetric case (i.e. $\alpha_1\le b \leq \alpha_1+\alpha_2$) is then (replacing $b$ by $\alpha_1+\alpha_2-b$ in the above rational function): 
\begin{equation}
h_b^{(\alpha_1,\alpha_2)}(p) = \frac{ p^{\alpha_1+\alpha_2-b+1} -1}{p-1}.
\end{equation}

The idea behind the method is the following. Let $F(x) = \sum a_nx^n$ and $G(x) = \sum b_nx^n$ be two convergent series. Then the series $H(x) = \sum a_nb_nx^n$  can be obtained from $F$ and $G$ as follows. Define the \textbf{convolution} of $F$ and $G$ denoted $F\star G (x)$ as the following contour integral around a small circle about  $0$ where $x$ is a complex number very small in magnitude.

\begin{align*}
F\star G (x) &= \frac{1}{2\pi i} \oint_{y\in B_{\epsilon}(0)} F(y)G(\frac{x}{y})\frac{dy}{y}\\
&= \frac{1}{2\pi i} \oint_{y\in B_{\epsilon}(0)} \sum a_nb_my^n\frac{x^m}{y^m}\frac{dy}{y} \\
&=\sum a_nb_mx^m \frac{1}{2\pi i} \oint_{y\in B_{\epsilon}(0)} y^{n-m-1}dy\\
&= \sum a_nb_nx^n.
\end{align*}
Thus $H(x) = F\star G (x).$ Now, we will prove equation (\ref{c1rk2}) using this method. Let $F_d(x,y)$ be the generating series
\begin{equation}
F_d(x,y) = \sum_{\substack{\alpha = (\alpha_1,\dots,\alpha_d) \\ \alpha_1\ge \cdots \alpha_d\ge 0  \\  0\le r\le \alpha_1+\dots+\alpha_d}} h_r^{(\alpha_1,\dots,\alpha_d)}(p)x^{\alpha}y^r 
\end{equation}
where $x = (x_1,\dots, x_d)$ and $x^\alpha = x_1^{\alpha_1}x_2^{\alpha_2}\cdots x_d^{\alpha_d}.$ Let's denote $h_r^{\alpha}(p)$ by $N_{\alpha}^{(d)}(r)$ or simply $N_{\alpha}(r)$ when the rank is understood. By \cite[Corollary]{TS}, we have the recurrence formula 
\begin{equation}\label{Stehrec}
    N_{\alpha}(r) = N_{\tilde{\alpha}}(r-1)+p^rN_{\hat{\alpha}}(r)
\end{equation} where $\hat{\alpha} = (\alpha_2,\dots,\alpha_d)$ and $\tilde{\alpha} = \alpha$ with $-1$ added to the $k^{th}$ position where 

$$k=\begin{cases}
  1  & \text{if } \alpha_1>\alpha_2\\    
  2  & \text{if } \alpha_1=\alpha_2>\alpha_3 \\    
  3  & \text{if } \alpha_1=\alpha_2=\alpha_3>\alpha_4\\
  \dots & \dots 
\end{cases}$$

Note $N_{\alpha}^{(1)}(r) =1$ if $0\le r\le \alpha.$ So, 
\begin{equation}
    F_1(x,y) = \frac{1}{(1-x)(1-xy)}.
\end{equation} For rank 2, let's break down $F_2(x_1,x_2,y)$ as follows: 
\begin{align}
\label{F2}
F_2(x_1,x_2,y) &=  \sum_{\substack{ \alpha_1 \ge \alpha_2\ge 0  \\  0\le r\le \alpha_1+\alpha_2}} N_{\alpha_1,\alpha_2}^{(2)}(r)x_1^{\alpha_1}x_2^{\alpha_2}y^r  
= F_2^{(0)}+F_2^{(1)}
\end{align} where 
\begin{align}
    F_2^{(0)}(x_1,x_2,y) &=  \sum_{\alpha_1 = \alpha_2\ge 0} N_{\alpha_1,\alpha_2}^{(2)}(r)x_1^{\alpha_1}x_2^{\alpha_2}y^r  \\
    F_2^{(1)}(x_1,x_2,y) &=  \sum_{\alpha_1 > \alpha_2\ge 0} N_{\alpha_1,\alpha_2}^{(2)}(r)x_1^{\alpha_1}x_2^{\alpha_2}y^r  
\end{align} Let's first compute $F_2^{(0)}(x_1,x_2,y). $
\begin{align*}
    F_2^{(0)}(x_1,x_2,y) &= \sum_{\alpha=0}^\infty \sum_r N_{\alpha,\alpha}^{(2)}(r)(x_1x_2)^{\alpha}y^r  \\
    &= \sum_{\alpha} \sum_r N_{\alpha,\alpha-1}^{(2)}(r-1)(x_1x_2)^{\alpha}y^r + \sum_{\alpha} p^rN_{\alpha}^{(1)}(r)(x_1x_2)^{\alpha}y^r\\
    &= \sum_{\alpha} \sum_r N_{\alpha-1,\alpha-1}^{(2)}(r-2)(x_1x_2)^{\alpha}y^r + \sum_{\alpha} \sum_r p^{r-1}N_{\alpha-1}^{(1)}(r-1)(x_1x_2)^{\alpha}y^r \\
    &+ \sum_{\alpha} p^rN_{\alpha}^{(1)}(r)(x_1x_2)^{\alpha}y^r\\
    &= \sum_{\alpha} \sum_r N_{\alpha,\alpha}^{(2)}(r)(x_1x_2)^{\alpha+1}y^{r+2} + \sum_{\alpha} \sum_r N_{\alpha}^{(1)}(r)(x_1x_2)^{\alpha+1}p^{r}y^{r+1}\\
    &+ \sum_{\alpha} N_{\alpha}^{(1)}(r)(x_1x_2)^{\alpha}p^ry^r\\
    &= x_1x_2y^2F_2^{(0)}(x_1,x_2,y) + x_1x_2yF_1(x_1x_2,py) + F_1(x_1x_2,py).
\end{align*} Solving for $F_2^{(0)},$ results in \begin{equation}
    F_2^{(0)}(x_1,x_2,y) = \frac{1+x_2x_2y}{(1-x_1x_2)(1-px_1x_2y)(1-x_1x_2y^2)}.
\end{equation} Similarly, we have 
\begin{align*}
    F_2^{(1)}(x_1,x_2,y) &= \sum_{\alpha_1>\alpha_2\ge 0} \sum_r N_{\alpha}^{(2)}(r)x_1^{\alpha_1}x_2^{\alpha_2}y^r  \\
    &=\sum_{\alpha_1>\alpha_2\ge 0} \left [ N_{\tilde{\alpha}}^{(2)}(r-1)+p^rN_{\hat{\alpha}}^{(1)}(r) \right ]x_1^{\alpha_1}x_2^{\alpha_2}y^r\\
    &=\sum_{\alpha_1>\alpha_2\ge 0} \left [ N_{\alpha_1-1,\alpha_2}^{(2)}(r-1)+p^rN_{\alpha_2}^{(1)}(r)  \right ]x_1^{\alpha_1}x_2^{\alpha_2}y^r\\
    &=\sum_{\alpha_1 \ge \alpha_2\ge 0}  N_{\alpha_1,\alpha_2}^{(2)}(r-1)x_1^{\alpha_1+1}x_2^{\alpha_2}y^{r+1} + \sum_{ \substack{ \alpha_1 > \alpha_2\ge 0\\ 0\le r\le \alpha_2}}  x_1^{\alpha_1}x_2^{\alpha_2}(py)^r\\
    &= x_1y\left [ F_2^{(1)}(x_1,x_2,y) + F_2^{(0)}(x_1,x_2,y)  \right ] + \sum_{ \alpha_1 > \alpha_2\ge 0}  x_1^{\alpha_1}x_2^{\alpha_2}\frac{1-(py)^{\alpha_2+1}}{1-py}.
\end{align*} Solving for $F_2^{(1)}(x_1,x_2,y),$ we have \begin{equation}
    F_2^{(1)}(x_1,x_2,y) = \frac{1+y-x_1y+x_1^2x_2y^2}{(1-x_1)(1-x_1y)(1-x_1x_2)(1-x_1x_2y)(1-qx_1x_2y)}.
\end{equation} Finally, combining the two sums, equation (\ref{F2}) becomes 

\begin{equation*}
F_2(x,y) = F_2(x_1,x_2,y) = \sum_{\substack{\alpha = (\alpha_1,\alpha_2) \\ \alpha_1\ge \alpha_2\ge 0  \\  0\le r\le \alpha_1+\alpha_2}} h_r^{(\alpha_1,\alpha_2)}(p)x_1^{\alpha_1}x_2^{\alpha_2}y^r\end{equation*}
\begin{equation}
= \frac{x_1^2x_2y^2 + x_1^2x_2y - x_1x_2y - 1}{(1-x_1)(1-x_1y)(1-x_1x_2)(1-x_1x_2y^2)(px_1x_2y-1)}.
\end{equation}

We will use the idea described before to pick out the rational function expression for case 1 from the general rational function expression of $F_2$ above. The final step will be extracting a general formula for the coefficients of the general term of the power series expansion of the resulting rational function using products of power series and change of variables.

Now consider small circles in the complex plane around 0 where the complex numbers $x_1,x_2,y$ are very small in magnitude.  We evaluate the following repeated contour integrals of the convolution of $F_2$ and $G$ by applying Cauchy's residue theorem to get the generating series for case 1.

Let 
\begin{equation}
G(x_1,x_2,y) = \sum_{ \alpha_1\ge \alpha_2\ge r\ge  0 } x_1^{\alpha_1}x_2^{\alpha_2}y^r = \frac{1}{(1-x_1)(1-x_1x_2)(1-x_1x_2y)}. 
\end{equation}
Then 
\begin{align*}
F_2\star G (x_1,x_2,y) &=\frac{1}{(2\pi i)^3} \oint_{u_3\in B_{\epsilon_3}(0)} \oint_{u_2\in B_{\epsilon_2}(0)} \oint_{u_1\in B_{\epsilon_1}(0)} F_2(u_1,u_2,v)G(\frac{x_1}{u_1},\frac{x_2}{u_2},\frac{y}{v})\frac{du_1}{u_1}\frac{du_2}{u_2}\frac{dv}{v}  \\  
&= \sum_{ \alpha_1\ge \alpha_2\ge r\ge  0 } h_r^{(\alpha_1,\alpha_2)}(p)x_1^{\alpha_1}x_2^{\alpha_2}y^r. 
\end{align*}

Using the rational functions for $F_2$ and $G$ and evaluating the integrals, we get that 
\begin{equation}
\sum_{ \alpha_1\ge \alpha_2\ge r\ge  0 } h_r^{(\alpha_1,\alpha_2)}(p)x_1^{\alpha_1}x_2^{\alpha_2}y^r = \frac{1}{(1-x_1)(1-x_1x_2)(1-x_1x_2y)(1-px_1x_2y)}.
\end{equation}
It remains to extract the coefficient of $x_1^{\alpha_1}x_2^{\alpha_2}y^r$ from the last rational function in order to determine $h_r^{(\alpha_1,\alpha_2)}(p)$ for the case $0\le r \le a_2.$ We do this as follows:
\begin{align*}
&\frac{1}{(1-x_1)(1-x_1x_2)(1-x_1x_2y)(1-px_1x_2y)}\\ 
&= \big( \sum_{k_1\ge 0} x_1^{k_1}\big)\big( \sum_{k_2\ge 0} (x_1x_2)^{k_2}\big)\big( \sum_{k_3\ge 0} (x_1x_2y)^{k_3}\big)\big( \sum_{k_4\ge 0} (px_1x_2y)^{k_3}\big)\\
&= \sum_{k_1,k_2,k_3,k_4\ge 0} p^{k_4}x_1^{k_1+k_2+k_3+k_4}x_2^{k_2+k_3+k_4}y^{k_3+k_4}.
\end{align*}

Let $r = k_3+k_4, a_2=k_2+k_3+k_4=k_2+r, a_1=k_1+k_2+k_3+k_4=k_1+a_2.$ The last sum above becomes 
\begin{align*}
\frac{1}{(1-x_1)(1-x_1x_2)(1-x_1x_2y)(1-px_1x_2y)} &=\sum_{a_1\ge a_2\ge r\ge 0} x_1^{a_1}x_2^{a_2}y^{r}\sum_{k_4=0}^r p^{k_4} \\
&= \sum_{a_1\ge a_2\ge r\ge 0} \frac{p^{r+1}-1}{p-1}x_1^{a_1}x_2^{a_2}y^{r}.
\end{align*}

Finally, by comparing coefficients, we get the required result:
$$h_r^{(\alpha_1,\alpha_2)}(p) = \frac{p^{r+1}-1}{p-1}.$$
\begin{remark}
In comparison, we can quickly prove the formulas for all cases of rank 2 using the recurrence used in the proof of Theorem \ref{rank3} as follows. First let's return to the notation of the previous sections, i.e. $a = (a_d, \dots,a_1)$ with $a_1\leq \dots \le a_d.$ 

Let $\lambda=(a_2, a_1)$ and $\lambda^{'}=( a_1).$\\
\text{Case 1:}  $0\le b \leq a_1$ 
$N_{b}(\lambda)=\sum_{i=0}^{b} p^{i}N_{i}(\lambda^{'})=\sum_{i=0}^{b} p^{i}=\frac{p^{b+1}-1}{p-1}.$\\
Note $N_{i}(\lambda^{'})=1$ or 0 according as $ i \leq a_1$ or $ i> a_1 ,$ respectively.\\
\text{Case 2:}  $a_1\le b \leq a_2$ 
$N_{b}(\lambda)=\sum_{i=0}^{a_1} p^{i}N_{i}(\lambda^{'})+\sum_{i=a_{1}+1}^{b} p^{i}N_{i}(\lambda^{'})=\sum_{i=0}^{a_1} p^{i}N_{i}(\lambda^{'})=\sum_{i=0}^{a_1} p^{i}=\frac{p^{a_1+1}-1}{p-1}.$\\
\text{Case 3:}  $a_2\le b \leq a_{1}+a_2$ (This case can also be deduced by replacing $b$ in case 1 by $a_1+a_2-b$)
$N_{b}(\lambda)=\sum_{i=0}^{a_1} p^{i}N_{i}(\lambda^{'})+\sum_{i=a_{1}+1}^{a_{2}} p^{i}N_{i}(\lambda^{'})+\sum_{i=a_{2}+1}^{b} p^{i}N_{i}(\lambda^{'})-\sum_{i=a_{1}+a_{2}+1-b}^{a_1}p^{i}N_{i}(\lambda^{'})$\\
$=\sum_{i=0}^{a_1} p^{i}N_{i}(\lambda^{'})-\sum_{i=a_{1}+a_{2}+1-b}^{a_1}p^{i}=\frac{p^{a_{1}+a_{2}+1-b}-1}{p-1}.$\\

In summary, we have
$$h_b^{(a_2,a_1)}(p)=\begin{cases}
  \frac{p^{b+1}-1}{p-1}  & \text{if } 0\le b \leq a_2\\    
  \frac{p^{a_2+1}-1}{p-1}  & \text{if } a_2\le b \leq a_1\\ 
  \frac{p^{a_1+a_2-b+1}-1}{p-1}  & \text{if } a_1\le b \leq a_1+a_2.
\end{cases}$$
\end{remark}

\textit{Acknowledgement.} The first named author would like to thank the City University of New York HPCC facilities for letting us verify every formula using the SageMath and Mathematica software on their HPC systems based at the College of Staten Island, New York City. 

\bibliographystyle{hplain}

\begin{thebibliography}{10}


\bibitem{GB96} Bhowmik, G., 1996. 
\newblock Evaluation of divisor functions of matrices. 
\newblock {\em Acta Arithmetica}, 74(2), pp.155-159.

\bibitem{GBJW} Bhowmik, G. and Wu, J., 1997.
\newblock On the asymptotic behaviour of the number of subgroups of finite abelian groups.
\newblock {\em Archiv der Mathematik}, 69(2), pp.95-104.

\bibitem{GBJW01} Bhowmik, G. and Wu, J., 2001. Zeta function of subgroups of abelian groups and average orders. Journal fur die Reine und Angewandte Mathematik, pp.1-16.


\bibitem{GB} Birkhoff, G., 1935. 
\newblock Subgroups of abelian groups. 
\newblock {\em Proceedings of the London Mathematical Society}, 2(1), pp.385-401.

\bibitem{butler}
Lynne~M. Butler.
\newblock A unimodality result in the enumeration of subgroups of a finite
  abelian group.
\newblock {\em Proc. Amer. Math. Soc.}, 101(4):771--775, 1987.

\bibitem{GC} Calugareanu, G., 2004. The total number of subgroups of a finite abelian group. {\em Scientiae Mathematicae Japonicae}, 60(1), pp.157-168.

\bibitem{CCL} C.Y. Chew, A.Y.M. Chin, and C.S. Lim. An explicit formula for the number of subgroups of a finite abelian p-group of rank 4. Unpublished

\bibitem{CKK} Chinta, G., Kaplan, N. and Koplewitz, S., 2017. The cotype zeta function of $\mathbb {Z}^ d$. arXiv preprint arXiv:1708.08547.

\bibitem{GLP} Goldfeld, D., Lubotzky, A. and Pyber, L., 2004. Counting congruence subgroups. {\em Acta Mathematica}, 193(1), pp.73-104.

\bibitem{HT} Hampejs, M. and Tóth, L., 2013. On the subgroups of finite Abelian groups of rank three. arXiv preprint arXiv:1304.2961.

\bibitem{HH} Hilton, H., 1907. On Sub‐Groups of a Finite Abelian Group. {\em Proceedings of the London Mathematical Society}, 2(1), pp.1-5.

\bibitem{YH} Hironaka, Y., 2017. Zeta functions of finite groups by enumerating subgroups. {\em Communications in Algebra}, 45(8), pp.3365-3376.

\bibitem{AI} Ivić, A., 1997. On the number of subgroups of finite abelian groups. {\em Journal de théorie des nombres de Bordeaux}, 9(2), pp.371-381.

\bibitem{IGM} Macdonald, I.G., 1998. Symmetric functions and Hall polynomials. {\em Oxford university press}.

\bibitem{MI04} Miller, G.A., 1904. On the subgroups of an abelian group. {\em The Annals of Mathematics}, 6(1), pp.1-6.
 
\bibitem{MI39} G. A. Miller, Number of the subgroups of any given abelian group. {\em Proc. Nat. Acad. Sci. U.S.A}. vol. 25 (1939) pp. 256-262.

\bibitem{OH} J.M. OH,  An explicit formula for the number of subgroups of a finite abelian p-group up to rank 3. {\em Commun. Math. Soc}. 28,649-667(2013).

\bibitem{SM} SageMath, the Sage Mathematics Software System (Version 8.1), The Sage Developers, 2018, http://www.sagemath.org.


\bibitem{SA} Amit Sehgal and Yogesh Kumar, On the number of subgroups of finite abelian $Z_m\otimes Z_n$. {\em  International Journal of Algebra}, 7 no. 19 (2013), 915-923.

\bibitem{PKASSS} Amit Sehgal , Sarita Sehgal and P.K. Sharma, The number of subgroups of a finite abelian p-group of rank 3. {\em Journal of Calcutta Mathematical Society}, Vol. 12, No. 2, 2016, 137-152.

\bibitem{RS} R. Schmidt, Subgroup Lattices of Groups. {\em de Gruyter Expositions in Mathematics 14, de Gruyter, Berlin}, (1994).

\bibitem{VNS} Shokuev, V.N., 1972. An expression for the number of subgroups of a given order of a finite p-group. Mathematical Notes, 12(5), pp.774-778.

\bibitem{TS} Stehling, T., 1992. On computing the number of subgroups of a finite abelian group. {\em Combinatorica}, 12(4), pp.475-479.

\bibitem{MT} T$\ddot{a}$rn$\ddot{a}$uceanu , M., 2010. An arithmetic method of counting the subgroups of a finite abelian group. {\em Bulletin mathématique de la Société des Sciences Mathématiques de Roumanie}, pp.373-386.

\bibitem{LT} T$\acute{o}$th, L., 2016. The number of subgroups of the group $\mathbb {Z} _m\times\mathbb {Z} _n\times\mathbb {Z} _r\times\mathbb {Z} _s$. arXiv preprint arXiv:1611.03302.


\end{thebibliography}

\end{document}